\documentclass[10pt,twoside]{siamltex}
\usepackage{amsmath}
\newcommand{\adj}{\text{adj}}
\newcommand{\nctwo}{\text{NC}^2}
\newtheorem{example}[theorem]{Example}
\title{Berkowitz's Algorithm and Clow Sequences}
\author{Michael Soltys\thanks{%
McMaster University,
Department of Computing and Software,
1280 Main Street West,
Hamilton, Ontario L8S 4K1, CANADA,
{\tt <soltys@mcmaster.ca>}.}}
\date{\today}
\begin{document}
\maketitle
\begin{abstract}
We present a combinatorial interpretation of Berkowitz's algorithm.
Berkowitz's algorithm, defined in~\cite{berk}, is the fastest known
parallel algorithm for computing the characteristic polynomial of a
matrix.  Our combinatorial interpretation is based on ``loop covers''
introduced by Valiant in \cite{valiant}, and ``clow sequences,''
defined in \cite{clows}.  Clow sequences turn out to capture very
succinctly the computations performed by Berkowitz's algorithm, which
otherwise is quite difficult to analyze.  The main contribution of
this paper is a proof of correctness of Berkowitz's algorithm in terms
of clow sequences (Theorem~\ref{thm:berk-clow}).
\end{abstract}

\section{Introduction}

Berkowitz's algorithm is the fastest known parallel algorithm for
computing the characteristic polynomial (char poly) of a matrix (and
hence for computing the determinant, the adjoint, and the inverse of a
matrix, if it exists).  It can be formalized with small boolean
circuits of depth $O(\log^2)$ in the size of the underlying matrix.
We shall describe precisely in the next section what we mean by
``small'' and ``depth,'' but the idea is that the circuits have
polynomially many gates in $n$, for an $n\times n$ matrix $A$, and the
longest path in these circuits is a constant multiple of $\log^2n$.

There are two other fast parallel algorithms for computing the
coefficients of the characteristic polynomial of a matrix: Chistov's
algorithm and Csanky's algorithm.  Chistov's algorithm is more
difficult to formalize, and Csanky's algorithm works only for fields
of char 0; see~\cite[section 13.4]{von} for all the details about
these two algorithms.

The author's original motivation for studying Berkowitz's algorithm was
the proof complexity of linear algebra.  Proof Complexity deals with
the complexity of formal mathematical derivations, and it has
applications in lower bounds and automated theorem proving.  In
particular, the author was interested in the complexity of derivations
of matrix identities such as $AB=I\rightarrow BA=I$ (right matrix
inverses are left inverses).   These identities have been proposed by
Cook as candidates for separating the Frege and Extended Frege proof
systems.   Proving (or disproving) this separation is one of the
outstanding problems in propositional proof complexity
(see~\cite{krajicek} for a comprehensive exposition of this area of
research).

Thus, we were interested in an algorithm that could compute inverses
of matrices, of the lowest complexity possible.  Berkowitz's algorithm
is ideal for our purposes for several reasons: 
\begin{itemize}
\item as was mentioned above, it has the lowest known complexity for
computing the char poly (we show in the next section that it can be
formalized with uniform $\nctwo$ circuits: small circuits of small
depth), 
\item it can be easily expressed with iterated matrix products, and
hence it lends itself to an easy formalization in first order logic
(with three sorts: indices, field elements, and matrices,
see~\cite{soltys}), 
\item and it is field independent, as the algorithm does not require
divisions, and hence Berkowitz's algorithm can compute char
polynomials over any commutative ring. 
\end{itemize}

Standard algorithms in linear algebra, such as Gaussian Elimination,
do not yield themselves to parallel computations.  Gaussian
Elimination is a sequential polynomial time algorithm, and hence it
falls in a complexity class far above the complexity of Berkowitz's
algorithm.  Furthermore, Gaussian Elimination requires divisions,
which are messy to formalize, and are not ``field independent''.  The
cofactor expansion requires computing $n!$ many terms, so it is an
exponential algorithm, and hence not tractable. 

From the beginning, we were interested in proving the correctness of
Berkowitz's algorithm within its own complexity class.  That is, for
our applications, we wanted to give a proof of correctness where the
computations were not outside the complexity class of Berkowitz's
algorithm, but rather inside $\nctwo$, meaning that the proof of
correctness should use iterated matrix products as its main engine for
computations.  This turned out to be a very difficult problem.

The original proof of correctness of Berkowitz's algorithm relies on
Samuelson's Identity (shown in the next section), which in turn relies
on Lagrange's Expansion, which is widely infeasible (as it requires
summing up $n!$ terms, for a matrix of size $n\times n$).  We managed
to give a feasible (polytime) proof of correctness in~\cite{soltys},
but the hope is that it is possible to give a proof of correctness
which does not need polytime concepts, but rather concepts from the
class $\nctwo$.   Note that by correctness, we mean that we can prove
the main properties of the char poly, which are: the Cayley-Hamilton
Theorem, and the multiplicativity of the determinant; all other
``universal'' properties follows directly from these two.

Hence our interest in understanding the workings of Berkowitz's
algorithm.  In this paper, we show that Berkowitz's algorithm computes
sums of the so called ``clow sequences.''  These are generalized
permutations, and they seem to be the conceptually cleanest way of
showing what is going on in Berkowitz's algorithm.  Since clow
sequences are generalized permutations, they do not lend themselves
directly to a feasible proof.  However, the hope is that by
understanding the complicated cancellations of terms that take place in
Berkowitz's algorithm,  we will be able to assert properties of
Berkowitz's algorithm which {\em do} have $\nctwo$ proofs, and which
imply the correctness of the algorithm.  Clow sequences expose very
concisely the cancellations of terms in Berkowitz's algorithm.

The main contribution of this paper is given by
Theorem~\ref{thm:berk-clow}, where we show that Berkowitz's algorithm
computes sums of clow sequences.  The first combinatorial
interpretation of Berkowitz's algorithm was given by Valiant
in~\cite{valiant}, and it was given in terms of ``loop covers,'' which
are similar to clow sequences.  This was more of an observation,
however, and not many details were given.  We give a detailed
inductive proof of correctness of Berkowitz's algorithm in terms of
clow sequences, introduced in~\cite{clows}. 

\section{Berkowitz's Algorithm}

Berkowitz's algorithm computes the coefficients of the char polynomial
of a matrix $A$, $p_A(x)=\det(xI-A)$, by computing iterated matrix
products, and hence it can be formalized in the complexity class
$\nctwo$.

The complexity class $\nctwo$ is the class of problems (parametrized
by $n$---here $n$ is the input size parameter) that can be computed
with uniform boolean circuits (with gates AND, OR, and NOT), of
polynomial size in $n$ (i.e., polynomially many gates in the input
size), and $O(\log^2n)$ depth (i.e., the longest path from an input
gate to the circuit to the output gate is in the order of $\log^2n$).  

For example, matrix powering is known to be in $\nctwo$.  The reason
is that the product of two matrices can be computed with boolean
circuits of polynomial size and logarithmic depth (i.e., in
$\text{NC}^1$), and the $n$-th power of a matrix can be obtained by
repeated squaring (squaring $\log n$ many times for a matrix of size
$n\times n$).  Cook defined\footnote{See (\cite{parallel}) for a
comprehensive exposition of the parallel classes $\text{NC}^i$, POW,
and related complexity classes.  This paper contains the details of
all the definitions outlined in the above paragraphs.} the complexity
class POW to be the class of problems reducible to matrix powering,
and showed that $\text{NC}^1\subseteq\text{POW}\subseteq\nctwo$.  Note
that every time we make the claim that Berkowitz's algorithm can be
formalized in the class $\nctwo$, we could be making a stronger claim
instead by saying that Berkowitz's algorithm can be formalized in the
class POW.

Berkowitz's algorithm computes the char polynomial of a matrix with
iterated matrix products.  Iterated matrix products can easily be
reduced to the problem of matrix powering: place the
$A_1,A_2,\ldots,A_n$ above the main diagonal of a new matrix $B$ which
is zero everywhere else, compute $B^n$, and extract $A_1A_2\cdots A_n$
from the upper-right corner block of $B^n$.   Hence, since Berkowitz's
algorithm can be computed with iterated matrix products (as we show in
Definition~\ref{dfn:berk} below), it follows that Berkowitz's
algorithm can be formalized in $\text{POW}\subseteq\nctwo$.  The
details are in Lemma~\ref{lem:nctwo} below.

The main idea in the standard proof of Berkowitz's algorithm
(see~\cite{berk}) is Samuelson's identity, which relates the char
polynomial of a matrix to the char polynomial of its principal
sub-matrix.  Thus, the coefficients of the char polynomial of an
$n\times n$ matrix $A$ below are computed in terms of the coefficients
of the char polynomial of $M$:
$$
A=\left(\begin{array}{lc}
a_{11} & R \\
S      & M
\end{array}\right)
$$
where $R,S$ and $M$ are $1\times (n-1)$, $(n-1)\times 1$ and
$(n-1)\times (n-1)$ sub-matrices, respectively.

\begin{lemma}[Samuelson's Identity]\label{lem:samuelson}
{\rm Let $p(x)$ and $q(x)$ be the char polynomials of $A$ and $M$,
respectively.  Then:
$$
p(x)=(x-a_{11})q(x)-R\cdot \adj(xI-M)\cdot S
$$}
\end{lemma}

Recall that the adjoint of a matrix $A$ is the transpose of the matrix
of cofactors of $A$; that is, the $(i,j)$-the entry of $\text{adj}(A)$
is given by $(-1)^{i+j}\det(A[j|i])$.  Also recall that $A[k|l]$ is
the matrix obtained from $A$ by deleting the $k$-th row and the $l$-th
column.  We also make up the following notation: $A[-|l]$ denotes that
only the $l$-th column has been deleted.  Similarly, $A[k|-]$ denotes
that only the $k$-th row has been deleted, and $A[-|-]=A$.

\begin{proof}(Lemma~\ref{lem:samuelson})
\begin{align*}
p(x) 	&= \det(xI-A) \\
	&= \det\left(\begin{array}{cc}
		x-a_{11} & -R \\
		-S       & xI-M
	   \end{array}\right) \\
\intertext{using the cofactor expansion along the first row:}
	&= (x-a_{11})\det(xI-M)+
	\sum_{j=1}^{n-1}(-1)^j(-r_j)\det(\underbrace{-S(xI-M)[-|j]}_{(*)})
\intertext{where $R=(r_1 r_2 \ldots r_{n-1})$, and the matrix
indicated by $(*)$ is given as follows: the first column is $S$, and
the remaining columns are given by $(xI-M)$ with the $j$-th column
deleted.  We expand $\det(-S(xI-M)[-|j])$ along the first column,
i.e., along the column $S=(s_1 s_2 \ldots s_{n-1})^T$ to obtain:}
	&= (x-a_{11})q(x)+
		\sum_{j=1}^{n-1}(-1)^j(-r_j)
		\sum_{i=1}^{n-1}(-1)^{i+1}(-s_i)\det(xI-M)[i|j] \\
\intertext{and rearranging:}
	&= (x-a_{11})q(x)-\sum_{i=1}^{n-1}\left(
		\sum_{j=1}^{n-1}r_j(-1)^{i+j}\det(xI-M)[i|j]
		\right)s_i \\
	&= (x-a_{11})q(x)-R\cdot \text{adj}(xI-M)\cdot S
\end{align*}
and we are done.
\end{proof}

\begin{lemma}\label{lem:adjoint}
{\rm Let $q(x)=q_{n-1}x^{n-1}+\cdots+q_1x+q_0$ be the char polynomial
of $M$, and let:
\begin{equation}\label{eq:adjoint}
B(x)=\sum_{k=2}^n(q_{n-1}M^{k-2}+\cdots+q_{n-k+1}I)x^{n-k}
\end{equation}
Then $B(x)=\text{adj}(xI-M)$.}
\end{lemma}

\begin{example}
{\rm If $n=4$, then 
\begin{align*}
B(x) &= Iq_3x^2+(Mq_3+Iq_2)x+(M^2q_3+Mq_2+Iq_1)
\end{align*}}
\end{example}

\begin{proof}(Lemma~\ref{lem:adjoint}) 
First note that:
$$
\text{adj}(xI-M)\cdot (xI-M)=\det(xI-M)I= q(x)I
$$
Now multiply $B(x)$ by $(xI-M)$, and using the Cayley-Hamilton
Theorem, we can conclude that $B(x)\cdot (xI-M)=q(x)I$.  Thus, the
result follows as $q(x)$ is not the zero polynomial; i.e., $(xI-M)$ is
{\em not} singular.
\end{proof}

From Lemma~\ref{lem:samuelson} and Lemma~\ref{lem:adjoint} we have
the following identity which is the basis for Berkowitz's algorithm:
\begin{equation}\label{eq:berk}
p(x)=(x-a_{11})q(x)-R\cdot B(x)\cdot S
\end{equation}

Using~(\ref{eq:berk}), we can express the char poly of a matrix as
iterated matrix product.  Again, suppose that $A$ is of the form:
$$
\left(\begin{array}{cc}
a_{11} & R \\
S      & M
\end{array}\right)
$$

\begin{definition}
{\rm We say that an $n\times m$ matrix is \emph{Toeplitz} if the
values on each diagonal are the same. We say that a matrix is
\emph{upper triangular} if all the values below the main diagonal are
zero.  A matrix is \emph{lower triangular} if all the values above the
main diagonal are zero.}
\end{definition}

If we express equation~(\ref{eq:berk}) in matrix form we obtain:
\begin{equation}\label{eq:berk-matrix}
p=C_1q
\end{equation}
where $C_1$ is an $(n+1)\times n$ Toeplitz lower triangular matrix,
and where the entries in the first column are defined as follows:
\begin{equation}\label{eq:def_of_C}
c_{i1}=\begin{cases}
1 & \text{if $i=1$} \\
-a_{11} & \text{if $i=2$} \\
-(RM^{i-3}S) & \text{if $i\geq 3$}
\end{cases}
\end{equation}

\begin{example}
{\rm If $A$ is a $4\times 4$ matrix, then $p=C_1q$ is given by:
$$
\left(\begin{array}{c} p_{4} \\ p_{3} \\ p_2 \\ p_1 \\ p_0
\end{array}\right)=
\left(\begin{array}{cccc}
1 & 0 & 0 & 0 \\
-a_{11} & 1 & 0 & 0 \\
-RS & -a_{11} & 1 & 0 \\
-RMS & -RS & -a_{11} & 1 \\
-RM^2S & -RMS & -RS & -a_{11} \\
\end{array}\right)
\left(\begin{array}{c} q_{3} \\ q_{2} \\ q_1 \\ q_0
\end{array}\right)
$$}
\end{example}

Berkowitz's algorithm consists in repeating this for $q$ (i.e., $q$
can itself be computed as $q=C_2r$, where $r$ is the char polynomial
of $M[1|1]$), and so on, and eventually expressing $p$ as a product of
matrices:
$$
p=C_1C_2\cdots C_n
$$
We provide the details in the next definition.

\begin{definition}[Berkowitz's algorithm]\label{dfn:berk}
{\rm Given an $n\times n$ matrix $A$, over any field $K$, {\em
Berkowitz's algorithms} computes an $(n+1)\times 1$ column vector
$p_A$ as follows: 

Let $C_j$ be an $(n+2-j)\times (n+1-j)$ Toeplitz and lower-triangular
matrix, where the entries in the first column are define as follows: 
\begin{equation}\label{eq:entries-general}
\begin{cases}
1          & \text{if $i=1$} \\
-a_{jj}    & \text{if $i=2$} \\
-R_jM_j^{i-3}S_j & \text{if $3\le i\le n+2-j$}
\end{cases}
\end{equation}
where $M_j$ is the $j$-th principal sub-matrix, so $M_1=A[1|1]$,
$M_2=M_1[1|1]$, and in general $M_{j+1}=M_j[1|1]$, and $R_j$ and $S_j$
are given by:
$$
\left(\begin{array}{cccc} a_{j(j+1)}&a_{j(j+2)}&\ldots&a_{jn}
\end{array}\right)\quad\text{and}\quad
\left(\begin{array}{cccc}
a_{(j+1)j}&a_{(j+2)j}&\ldots&a_{nj} \end{array}\right)^t
$$
respectively.   Then:
\begin{equation}\label{eq:product}
p_A=C_1C_2\cdots C_n
\end{equation}}
\end{definition}
Note that Berkowitz's algorithm is field independent (there are no
divisions in the computation of $p_A$), and so all our results are
field independent.

\begin{lemma}\label{lem:nctwo}
{\rm Berkowitz's algorithm is an $\nctwo$ algorithm.}
\end{lemma}

\begin{proof}
This follows from~(\ref{eq:product}): $p_A$ is given as a product of
matrices, each $C_i$ can be computed independently of the other
$C_j$'s, so we have a sequence of $C_1,C_2,\ldots,C_n$ matrices,
independently computed, so we can compute their product with repeated
squaring of the matrix $B$, which is constructed by placing the
$C_i$'s above the main diagonal of an otherwise all zero matrix.

Now the entries of each $C_i$ can also be computed using matrix
products, again independently of each other.  In fact, we can compute
the $(i,j)$-th entry of the $k$-th matrix very quickly as
in~(\ref{eq:entries-general}).  

Finally, we can compute additions, additive inverses, and products of
the underlying field elements (in fact, more generally, of the
elements in the underlying commutative ring, as we do not need
divisions in this algorithm).  We claim that these operations can be
done with small $\text{NC}^1$ circuits (this is certainly true for the
standard examples: finite fields, rationals, integers, etc.).  

Thus we have ``three layers'': one layer of $\text{NC}^1$ circuits,
and two layers of $\nctwo$ circuits (one layer for computing the
entries of the $C_j$'s, and another layer for computing the product of
the $C_j$'s), and so we have (very uniform) $\nctwo$ circuits that
compute the a column vector with the coefficients of the char
polynomial of a given matrix.
\end{proof}

In this section we showed that Berkowitz's algorithm computes the
coefficients of the char polynomial correctly, by first proving
Samuelson's Identity, and then using the Cayley-Hamilton Theorem, to
finally obtain that equation~(\ref{eq:product}) computes the
coefficients of the char polynomial correctly.  This is a very
indirect approach, and we loose insight into what is actually being
computed when we are presented with equation~(\ref{eq:product}).
However, the underlying fact is the Lagrange expansion of the
determinant (that's how Samuelson's Identity, and the Cayley-Hamilton
Theorem are proved).  In the next section, we take
equation~(\ref{eq:product}) naively, and we give a combinatorial proof
of its correctness with clow sequences and the Lagrange expansion.

\section{Clow Sequences}

First of all, a ``clow'' is an acronym for ``closed walk.''  Clow
sequences (introduced in \cite{clows}, based on ideas in
\cite{straubing}), can be thought of as generalized permutations.
They provide a very good insight into what is actually being computed
in Berkowitz's algorithm.  

In the last section, we derived Berkowitz's algorithm from Samuelson's
Identity and the Cayley-Hamilton Theorem.  However, both these
principles are in turn proved using Lagrange's expansion for the
determinant.  Thus, this proof of correctness of Berkowitz's algorithm
is indirect, and it does not really show what is being computed in
order to obtain the char polynomial.

To see what is being computed in Berkowitz's algorithm, and to
understand the subtle cancellations of terms, it is useful to look at
the coefficients of the char polynomial of the determinant of a matrix
$A$ as given by determinants of minors of $A$.  To define this notion
precisely, let $A$ be an $n\times n$ matrix, and define
$A[i_1,\ldots,i_k]$, where $1\le i_1<i_2<\cdots<i_k\le n$, to be the
matrix obtained from $A$ by deleting the rows and columns numbered by
$i_1,i_2,\ldots,i_k$.  Thus, using this notation, $A[1|1]=A[1]$, and
$A[2,3,8]$ would be the matrix obtained from $A$ by deleting rows and
columns $2,3,8$.

Now, it is not difficult to show from the Lagrange's expansion of
$\det(xI-A)$ that $p_{n-k}$, where $p_n,p_{n-1},\ldots,p_0$ are the
coefficients of the char polynomial of $A$, is given by the following
formula:
\begin{equation}\label{eq:4}
p_k=\sum_{1\le i_1<i_2<\cdots<i_k\le n}
\det(A[i_1,i_2,\ldots,i_k])
\end{equation}
Since $\det(A[i_1,i_2,\ldots,i_k])$ can be computed using the
Lagrange's expansion, it follows from~(\ref{eq:4}), that each
coefficient of the char polynomial can be computed by summing over
permutations of minors of $A$:
\begin{equation}\label{eq:5}
p_{n-k}=\sum_{1\le i_1<i_2<\cdots<i_{n-k}\le n}
\sum_{\sigma\in\{j_1,j_2,\ldots,j_k\}}\text{sign}(\sigma)
a_{j_1\sigma(j_1)}a_{j_2\sigma(j_2)}\cdots a_{j_k\sigma(j_k)}
\end{equation}
Note that the set $\{j_1,j_2,\ldots,j_k\}$ is understood to be the
complement of the set $\{i_1,i_2,\ldots,i_{n-k}\}$ in
$\{1,2,\ldots,n\}$, and $\sigma$ is a permutation in $S_{k}$, that
permutes the elements in the set $\{j_1,j_2,\ldots,j_k\}$.  We
re-arranged the subscripts in~(\ref{eq:5}) to make the expression more
compatible later with clow sequences.  Note that when $k=n$, we are
simply computing the determinant, since in that case the first
summation is empty, and the second summation spans over all the
permutations in $S_n$:
$$
\det(A)=\sum_{\sigma\in S_n}\text{sign}(\sigma)a_{1\sigma(1)}\cdots
a_{n\sigma(n)}
$$
Finally, note that if $k=0$, then the second summation is empty, and
there is just one sequence that satisfies the condition of the first
summation (namely $1<2<\cdots<n$), so the result is 1 by convention.

We can interpret $\sigma\in S_n$ as a directed graph $G_\sigma$ on $n$
vertices: if $\sigma(i)=j$, then $(i,j)$ is an edge in $G_\sigma$, and
if $\sigma(i)=i$, then $G_\sigma$ has the self-loop $(i,i)$.  

\begin{example}\label{ex:1}
{\rm The permutation given by:
$$
\sigma=\left(\begin{array}{cccccc}
1 & 2 & 3 & 4 & 5 & 6 \\
3 & 1 & 2 & 4 & 6 & 5
\end{array}\right)
$$
corresponds to the directed graph $G_\sigma$ with 6 nodes and the
following edges: 
$$
\{(1,3),(2,1),(3,2),(4,4),(5,6),(6,5)\}
$$
where $(4,4)$ is a self-loop.}
\end{example}

Given a matrix $A$, define the weight of $G_\sigma$, $w(G_\sigma)$, as
the product of $a_{ij}$'s such that $(i,j)\in G_\sigma$.  So
$G_\sigma$ in example~\ref{ex:1} has a weight given by:
$w(G_\sigma)=a_{13}a_{21}a_{32}a_{44}a_{56}a_{65}$.  Using the new
terminology, we can restate equation~(\ref{eq:5}) as follows:
\begin{equation}\label{eq:6}
p_{n-k}=\sum_{1\le i_1<i_2<\cdots<i_{n-k}\le n}
\sum_{\sigma\in\{j_1,j_2,\ldots,j_k\}}\text{sign}(\sigma)
w(G_\sigma)\tag{$\ref{eq:5}'$}
\end{equation}

The graph-theoretic interpretation of permutations gets us closer to
clow sequences.  The problem that we have with~(\ref{eq:6}) is that
there are too many permutations, and there is no (known) way of
grouping or factoring them, in such a way so that we can save
computing all the terms $w(\sigma)$, or at least so that these terms
cancel each other out as we go.  

The way to get around this problem is by generalizing the notion of
permutation (or {\em cycle cover}, as permutations are called in the
context of the graph-theoretic interpretation of $\sigma$).  Instead
of summing over cycle covers, we sum over clow sequences; the paradox
is that there are many more clow sequences than cycle covers, {\em
but} we can efficiently compute the sums of clow sequences (with
Berkowitz's algorithm), making a clever use of cancellations of terms
as we go along.  We now introduce all the necessary definitions,
following \cite{clows}.

\begin{definition}
{\rm A \emph{clow} is a walk $(w_1,\ldots,w_l)$ starting from vertex
$w_1$ and ending at the same vertex, where any $(w_i,w_{i+1})$ is an
edge in the graph.  Vertex $w_1$ is the least-numbered vertex in the
clow, and it is called the \emph{head} of the clow.  We also require
that the head occur only once in the clow.  This means that there is
exactly one incoming edge $(w_l,w_1)$, and one outgoing edge
$(w_1,w_2)$ at $w_1$, and $w_i\neq w_1$ for $i\neq 1$.  The {\em
length} of a clow $(w_1,\ldots,w_l)$ is~$l$.  Note that clows are not
allowed to be empty, since they always must have a head.}
\end{definition}

\begin{example}\label{exm:clow}
{\rm Consider the clow $C$ given by $(1,2,3,2,3)$ on four vertices.
The head of clow $C$ is vertex 1, and the length of $C$ is 6.} 
\end{example}

\begin{definition}
{\rm A \emph{clow sequence} is a sequence of clows $(C_1,\ldots,C_k)$,
where $\text{head}(C_1)<\ldots<\text{head}(C_k)$.  The {\em length} of
a clow sequence is the sum of the lengths of the clows (i.e., the
total number of edges, counting multiplicities).  Note that a cycle
cover is a special type of a clow sequence.}
\end{definition}


\begin{definition}\label{dfn:sign}
{\rm We define the {\em sign} of a clow sequence to be $(-1)^k$ where
$k$ is the number of clows in the sequence.}
\end{definition}

\begin{example}\label{ex:cancel}
{\rm We list the clow sequences associated with the three vertices
$\{1,2,3\}$.  We give the sign of the corresponding clow sequences in
the right-most column:
$$
\begin{array}{lll}
1.  & (1),(2),(3) & (-1)^{3+3}=1   \\
2.  & (1,2),(3)   & (-1)^{3+2}=-1  \\
3.  & (1,2,2)     & (-1)^{3+1}=1   \\
4.  & (1,2),(2)   & (-1)^{3+2}=-1  \\
5.  & (1),(2,3)   & (-1)^{3+2}=-1  \\
6.  & (1,2,3)     & (-1)^{3+1}=1   \\
7.  & (1,3,3)     & (-1)^{3+1}=1   \\
8.  & (1,3),(3)   & (-1)^{3+2}=-1  \\
9.  & (1,3,2)     & (-1)^{3+1}=1   \\
10. & (1,3),(2)   & (-1)^{3+2}=-1  \\
11. & (2,3,3)     & (-1)^{3+1}=1   \\
12. & (2,3),(3)   & (-1)^{3+2}=-1 
\end{array}
$$
Note that the number of permutations on $3$ vertices is $3!=6$, and
indeed, the clow sequences $\{3,4,7,8,11,12\}$ do {\em not} correspond
to cycle covers.   We listed these clow sequences which do not
correspond to cycle covers by pairs: $\{3,4\},\{7,8\},\{11,12\}$.
Consider the first pair: $\{3,4\}$.  We will later define the weight
of a clow (simply the product of the labels of the edges), but notice
that clow sequence 3 corresponds to $a_{12}a_{22}a_{21}$ and clow
sequence 4 corresponds to $a_{12}a_{21}a_{22}$, which is the same
value; however, they have opposite signs, so they cancel each other
out.  Same for pairs $\{7,8\}$ and $\{11,12\}$.  We make this informal
observation precise with the following definitions, and in
Theorem~\ref{thm:clow} we show that clow sequences which do not
correspond to cycle covers cancel out.}
\end{example}

Given a matrix $A$, we associate a weight with a clow sequence that is
consistent with the contribution of a cycle cover.  Note that we can
talk about clows and clow sequences independently of a matrix, but
once we associate weights with clows, we have to specify the
underlying matrix, in order to label the edges.  Thus, to make things
more precise, we will sometimes say ``clow sequences on $A$'' to
emphasize that the weights come from $A$.

\begin{definition}
{\rm Given a matrix $A$, the {\em weight of a clow} $C$, denoted
$w(C)$, is the product of the weights of the edges in the clow, where
edge $(i,j)$ has weight $a_{ij}$.}  
\end{definition}  

\begin{example}
{\rm Given a matrix $A$, the weight of clow $C$ in
example~\ref{exm:clow} is given by:
$$
w((1,2,3,2,3))=a_{12}a_{23}^2a_{32}a_{31}
$$}
\end{example}

\begin{definition}\label{dfn:weight} 
{\rm Given a matrix $A$, the {\em weight of a clow sequence} $C$,
denoted $w(C)$, is the product of the weights of the clows in $C$.
Thus, if $C=(C_1,\ldots,C_k)$, then: 
$$
w(C)=\prod_{i=1}^kw(C_i).
$$
We make the convention that an empty clow sequence has weight 1.
Since a clow must consist of at least one vertex, a clow sequence is
empty iff it has length zero.  Thus, equivalently, a clow sequence of
length zero has weight 1.  These statements will be important when we
link clow sequences with Berkowitz's algorithm.}
\end{definition}

\begin{theorem}\label{thm:clow}
{\rm Let $A$ be an $n\times n$ matrix, and let
$p_n,p_{n-1},\ldots,p_0$ be the coefficients of the char polynomial of
$A$ given by $\det(xI-A)$.  Then:
\begin{equation}\label{eq:summation}
p_{n-k}=\sum_{\mathcal{C}_k}\text{sign}(C)w(C)
\end{equation}
Where $\mathcal{C}_k=\{C|\text{$C$ is a clow sequence on $A$ of length
$k$}\}$.} 
\end{theorem}

\begin{proof}
We generalize the proof given in~\cite[pp.~5--8]{clows} for the case
$k=n$.  The main idea in the proof is that clow sequences which are
not cycle covers cancel out, just as in example~\ref{ex:cancel}, so
the contribution of clow sequences which are not cycles covers is zero.  

Suppose that $(C_1,\ldots,C_j)$ is a clow sequence in $A$ of length
$k$.  Choose the smallest $i$ such that $(C_{i+1},\ldots,C_j)$ is a
set of disjoint cycles.  If $i=0$, $(C_1,\ldots,C_j)$ is a cycle
cover.  Otherwise, if $i>0$, we have a clow sequence which is not a
cycle cover, so we show how to find another clow sequence (which is
also not a cycle cover) of the same weight and length, but opposite
sign.  The contribution of this pair to the summation
in~(\ref{eq:summation}) will be zero.

So suppose that $i>0$, and traverse $C_i$ starting from the head until
one of two possibilities happens: (i) we hit a vertex that is in
$(C_{i+1},\ldots,C_j)$, or (ii) we hit a vertex that completes a
simple cycle in $C_i$.  Denote this vertex by $v$.  In case (i), let
$C_p$ be the intersected clow ($p\ge i+1$), join $C_i$ and $C_p$ at
$v$ (so we merge $C_i$ and $C_p$).  In case (ii), let $C$ be the
simple cycle containing $v$: detach it from $C_i$ to get a new clow.

In either case, we created a new clow sequence, of {\em opposite sign
but same weight and same length $k$}.  Furthermore, the new clow
sequence is still not a cycle cover, and if we would apply the above
procedure to the new clow sequence, we would get back the original
clow sequence (hence our procedure defines an {\em involution} on the
set of clow sequences).
%
\end{proof}

In~\cite{valiant} Valiant points out that Berkowitz's algorithm
computes sums of what he calls ``loop covers.''  We show that
Berkowitz's algorithm computes sums of slightly restricted clow
sequences, which are nevertheless equal to the sums of all clow
sequences, and therefore, by Theorem~\ref{thm:clow}, Berkowitz's
algorithm computes the coefficients $p_{n-k}$ of the char polynomial
of $A$ correctly.  We formalize this argument in the next theorem,
which is the central result of this paper.

%
%
%
%

\begin{theorem}\label{thm:berk-clow}
{\rm Let $A$ be an $n\times n$ matrix, and let: 
$$
p_A=\left(\begin{array}{c}
p_n\\ p_{n-1}\\\vdots\\ p_0
\end{array}\right)
$$ 
as in defined by equation~\ref{eq:product};  that is, $p_A$ is the
result of running Berkowitz's algorithm on $A$.  Then, for $0\le i\le
n$, we have:
\begin{equation}\label{eq:3}
p_{n-i}=\sum_{\mathcal{C}_i}\text{sign}(C)w(C)
\end{equation}
where $\mathcal{C}_i=\{C|\text{$C$ is a clow sequence on $A$ of length
$i$}\}$.} 
\end{theorem}

%
%

Before we prove this theorem, we give an example.
\begin{example}\label{exm:last}
{\rm Suppose that $A$ is a $3\times 3$ matrix, $M=A[1|1]$ as usual,
and $p_3,p_2,p_1,p_0$ are the coefficients of the char poly of $A$ and
$q_2,q_1,q_0$ are the coefficients of the char poly or $M$, computed
by Berkowitz's algorithm.  Thus:
\begin{equation}\label{eq:B1}
\begin{split}
\left(\begin{array}{c}
p_3 \\ p_2 \\ p_1 \\ p_0
\end{array}\right)&=
\left(\begin{array}{ccc}
1 & 0 & 0 \\ -a_{11} & 1 & 0 \\ -RS & -a_{11} & 1 \\
-RMS & -RS & -a_{11}
\end{array}\right)
\left(\begin{array}{c}
q_2 \\ q_1 \\ q_0
\end{array}\right) \\
&=\left(\begin{array}{c}
q_2 \\ -a_{11}q_2+q_1 \\ -RSq_2-a_{11}q_1+q_0 \\
-RMSq_2-RSq_1-a_{11}q_0\quad(\ast)
\end{array}\right)
\end{split}
\end{equation}
We assume that the coefficients $q_2,q_1,q_0$ are given by sums of
clow sequences on $M$, that is, by clow sequences on vertices
$\{2,3\}$.  Using this assumption and equation~(\ref{eq:B1}), we show
that $p_3,p_2,p_1,p_0$ are given by clow sequences on $A$, just as in
the statement of Theorem~\ref{thm:berk-clow}. 

Since $q_2=1$, $p_3=1$ as well.  Note that $q_2=1$ is consistent with
our statement that it is the sum of restricted clow sequences of
length zero, since there is only one empty clow sequence, and its
weight is by convention 1 (see Definition~\ref{dfn:weight}).  

Consider $p_2$, which by definition is supposed to be the sum of clow
sequences of length one on all three vertices.  This is the sum of
clow sequences of length one on vertices 2 and 3 (i.e., $q_1$), plus
the clow sequence consisting of a single self-loop on vertex 1 with
weight $a_{11}$ and sign $(-1)^{1}=-1$ (see
Definition~(\ref{dfn:sign})).  Hence, the sum is indeed
$-a_{11}q_2+q_1$, as in equation~(\ref{eq:B1}) (again, $q_2=1$).

Consider $p_1$.  Since $p_1=p_{3-2}$, $p_1$ is the sum of clow
sequences of length two.  We are going to show that the term
$-RSq_2-a_{11}q_1+q_0$ is equal to the sum of clow sequences of length
2 on $A$.   First note that there is just one clow of length two on
vertices 2 and 3, and it is given by $q_0$.
There are two clows of length two which include a self loop at vertex
1. These clows correspond to the term $-a_{11}q_1$.  Note that the
negative sign comes from the fact that $q_1$ has a negative value, but
there are two clows per sequence, so the parity is even, according to
Definition~(\ref{dfn:sign}).
Finally, we consider the clow sequences of length two, where there is
no self loop at vertex 1.  Since vertex 1 must be included, there are
only two possibilities; these clows correspond to the term $-RSq_2$
which is equal to:
$$
-\left(\begin{array}{cc} a_{12} & a_{13} \end{array}\right)
\left(\begin{array}{c} a_{21} \\ a_{31} \end{array}\right)=
-a_{12}a_{21}-a_{13}a_{31}
$$
since $q_2=1$.

For $p_0$, the reader can add up all the clows by following
Example~(\ref{ex:cancel}).  One thing to notice, when tracing this
case, is that the summation indicated by $(\ast)$ includes only those
clow sequences which start at vertex 1.  This is because, the bottom
entry in equation~(\ref{eq:B1}), unlike the other entries, does not
have a 1 in the last column, and hence there is not coefficient from
the char poly of $M$ appearing by itself.  This is not a problem for
the following reason: if vertex 1 is not included in a clow sequence
computing the last entry, then that clow sequence will cancel out
anyways, since a clow sequence of length 3 that avoids the first
vertex, cannot be a cycle cover!  This observation will be made more
explicit in the proof below.}
\end{example}

\begin{proof}(Theorem~\ref{thm:berk-clow})
We prove this theorem by induction on the size of matrices.  The Basis
Case is easy, since if $A$ is a $1\times 1$ matrix, then $A=(a)$, so
$p_A=\left(\begin{array}{cc} 1 & -a \end{array}\right)$, so $p_1=1$,
and $p_0=-a$ which is $(-1)\times$ the sum of clow sequences of length
1.

In the Induction Step, suppose that $A$ is an $(n+1)\times(n+1)$
matrix and:
\begin{equation}\label{eq:pCq}
\left(\begin{array}{c} 
p_{n+1} \\ 
p_{n} \\ 
p_{n-1} \\
p_{n-2} \\
p_{n-3} \\
\vdots \\ 
p_0
\end{array}\right)\!\!=\!\!
\left(\begin{array}{cccc}
1 & 0 & 0 &  \ldots  \\
-a_{11} & 1 & 0 & \ldots \\
-RS & -a_{11} & 1 & \ldots \\
-RMS & -RS & -a_{11} & \ldots\\
-RM^2S & -RMS & -RS & \ldots \\
\vdots &\vdots &\vdots &\ddots  \\
-RM^{n-1}S & -RM^{n-2}S & -RM^{n-3}S & \ldots
\end{array}\right)
\left(\begin{array}{c} 
q_{n} \\ 
q_{n-1} \\ 
q_{n-2} \\ 
q_{n-3} \\ 
\vdots \\
q_0
\end{array}\right)
\end{equation}
By the Induction Hypothesis, $q_M=\left(\begin{array}{cccc} q_n &
q_{n-1} & \ldots & q_0 \end{array}\right)$ satisfies the statement of
the theorem for $M=A[1|1]$, that is, $q_{n-i}$ is equal to the sum of
clow sequences of length $i$ on $M=A[1|1]$.

Since $p_{n+1}=q_n$, $p_{n+1}=1$.  Since
$p_{n}=-a_{11}q_n+q_{n-1}=-a_{11}+q_{n-1}$ (as $q_n=1$), using the
fact that $q_{n-1}=$ the sum of clow sequences of length 1 on $M$, it
follows that $p_n=$ the sum of clow sequences of length 1 on $A$.

Now we prove this for general $n+1>i>1$, that is, we prove that
$p_{n+1-i}$ is the sum of clow sequences of length $i$ on $A$.  Note
that:
\begin{equation}\label{eq:pi}
p_{n+1-i}=-RM^{i-2}Sq_n-RM^{i-3}Sq_{n-1}-\cdots-
RSq_{n+2-i}-a_{11}q_{n+1-i}+q_{n-i}
\end{equation}
as can be seen by inspection from equation~(\ref{eq:pCq}).  Observe
that the $(i,j)$-th entry of $M^k$ is the sum of clows in $M$ that
start at vertex $i$ and end at vertex $j$ of length $k$, and
therefore, $RM^kS$ is the sum of clows in $A$ that start at vertex 1
(and of course end at vertex 1, and vertex 1 is never visited
otherwise), of length $k+2$.

Therefore, $RM^{i-2-j}Sq_{n-j}$, for $j=0,\ldots,i-2$, is the product
of the sum of clows of length $i-j$ (that start and end at vertex 1)
and the sum of clow sequences of length $j$ on $M$ (by the Induction
Hypothesis), which is just the sum of clow sequences of length $i$
where the first clow starts and ends at vertex 1, and has length
$i-j$.  Each clow sequence of length $i$ on $A$ starts off with a clow
anchored at the first vertex, and the second to last term of
equation~(\ref{eq:pi}), $-a_{11}q_{n+1-i}$, corresponds to the case
where the first clow is just a self loop.  Finally, the last term
given by $q_{n-i}$ contributes the clow sequences of length $i$ which
do {\em not} include the first vertex.

The last case is when $i=n+1$, so $p_0$, which is the determinant of
$A$, by Theorem~\ref{thm:clow}.  As was mentioned at the end of
Example~\ref{exm:last}, this is a {\em special} sum of clow sequences,
because the head of the first clow is always vertex 1.  Here is when
we invoke the proof of the Theorem~\ref{thm:clow}: the last entry,
$p_0$ can be shown to be the sum of clow sequences, where the head of
the first clow is always vertex 1, by following an argument analogous
to the one in the above paragraph.  However, this sum is still equal
in value to the sum of all clow sequences (of length $n+1$).  This is
because, if we consider clow sequences of length $n+1$, and there are
$n+1$ vertices, and we get a clow sequences $C$ which avoids the first
vertex, then we know that $C$ {\em cannot} be a cycle cover, and
therefore it will cancel out in the summation anyways, just as it was
shown to happen in the proof of Theorem~\ref{thm:clow}. 
\end{proof}

\section{Acknowledgments}

I would like to thank Stephen Cook, who introduced me to Berkowitz's
algorithm, and who supervised my PhD thesis, where we designed logical
theories for reasoning about Linear Algebra, and where we gave the
first feasible proof of the Cayley-Hamilton Theorem (and other
fundamental theorems of Linear Algebra) based on Berkowitz's
algorithm; see~\cite{soltys}.

\bibliographystyle{plain}
\bibliography{berk}

\end{document}